\documentclass{article}
\usepackage{amssymb,amsmath,eucal}

\def\R {{\Bbb R }}
\def\C {{\Bbb C }}
\def\K{{\Bbb K}}
\def\H{{\Bbb H}}
\def\Z{{\Bbb Z}}

\def\la{\langle}
\def\ra{\rangle}

\def\RA{\longrightarrow}

\def\const{{\rm const}}
\def\B{{\rm B}}

\def\OO{{\rm O}}

\def\SL{{\rm SL}}

\def\Sp{{\rm Sp}}

\def\U{{\rm U}}

\def\cP{\mathcal P}
\def\cM{\mathcal M}
\def\cD{\mathcal D}

\def\ov{\overline}
\def\phi{\varphi}
\def\epsilon{\varepsilon}
\def\kappa{\varkappa}
\def\le{\leqslant}
\def\ge{\geqslant}

\def\densitybc{\Bigl|\frac
     {\Gamma(b+is)\Gamma(c+is)}{\Gamma(2is)}\Bigr|^2}

\def\densityx{x^{b+c-1}(1+x)^{b-c}}

\def\wh{\widehat}

\def\F{\,\,{}_2F_1}
\def\FF{\,\,{}_3F_2}
\def\FFF{\,\,{}_4F_3}

\def\FBC{\F(b+is, b-is; b+c; -x)}

\def\intt{\int_0^\infty}

\renewcommand{\Re}{\mathop{\rm Re}\nolimits}
\renewcommand{\Im}{\mathop{\rm Im}\nolimits}

\newcounter{sec}

\newcounter{punct}[sec]

\def\punct{\refstepcounter{punct}{\arabic{sec}.\arabic{punct}.  }}

\newtheorem{theorem}{Theorem}[sec]
\newtheorem{proposition}[theorem]{Proposition}

\newtheorem{lemma}[theorem]{Lemma}

\newtheorem{corollary}[theorem]{Corollary}

\def\COUNTERS{\addtocounter{sec}{1}
              \setcounter{punct}{0}
          \setcounter{equation}{0}
          \setcounter{theorem}{0}
          \setcounter{problem}{0}
          }

\begin{document}

\begin{center}

{\bf\Large
Index hypergeometric integral transform}

\bigskip

\sc \large
Yury A.Neretin\footnote{Partially supported by the grant FWF, project P22122}
\bf\Large
  \end{center}

{\small This is a brief overview  of the index hypergeometric  transform 
(other terms for this integral operator are: Olevskii transform, Jacobi transform,
generalized Mehler--Fock transform). We discuss applications of this transform
to special functions and harmonic analysis.
The text is an addendum to the Russian edition of the book by G.~E.~Andrews, R.~Askey, 
and R.~Roy, Special Functions, Encycl. of Math. Appl. 71, Cambridge Univ. Press, 1999.}

\bigskip

As is well-known, the continuous analog of Fourier series is the Fourier transform.
It turns out that  expansions in the Jacobi polynomials also have
a continuous analog, namely the integral 'Jacobi transform'
(the terminology is discussed below).  The theory of this transform
is rich, but there is no detailed modern exposition
of this topic in the existing literature, this text also
has not such a purpose. 
Numerous additional facts are contained in works
  \cite{Che}, \cite{DS}, \cite{FK},
 \cite{Koo1}--\cite{Koo3},
\cite{Ner-index}--\cite{Ner-wilson},  \cite{YaL}.

There is a well-known Askey--Wilson hierarchy of hypergeometric orthogonal
polynomials (see \cite{KLS}). There is a parallel hierarchy
of hypegeometric integral transforms, see
 \cite{Che2},
 \cite{Gro}, \cite{DS}, \cite{Koe}, \cite{Mol},
 \cite{Ner-hahn}--\cite{Ner-kus}, \cite{Wimp},
\cite{YaL}, on multidimensional analogs, see
\cite{HO}, \cite{Mac}, \cite{Che}. 
Our topic is neither the simplest object of this hierarchy
(certainly the Hankel transform and the Kontorovich--Lebedev
transform are simpler), nor the most complicated
(see, for instance \cite{Gro}).
It is sufficiently simple to be versatile tool
for special functions (see below Section 2),
on the other hand it controls harmonic analysis
on hyperbolic symmetric spaces (i.e., 
Lobachevsky spaces and their complex and quaternionian analogs),
we briefly discuss this in Section 4, for more details see
 \cite{Koo2},
\cite{Ner-index}.

\section{The index hypergeometric transform}

\COUNTERS

\smallskip

{\bf\punct The Jacobi polynomials.} Consider the Jacobi
orthogonal system on the segment
 $[0,1]$,
$$
\cP_n^{\alpha,\beta}(x)=
\frac{(-1)^n\Gamma(n+\beta+1)}{\Gamma(\beta+1) n!}
\,\,
\F\left[
\begin{matrix}-n,n+\alpha+\beta+1\\
 \beta+1
 \end{matrix};
  x\right]
.$$
We have
\begin{multline}
\gamma_n:=
\|P_n^{\alpha,\beta}\|^2
= \\ =
\int_0^1 \cP_n^{\alpha,\beta}(x)^2\, x^\beta(1-x)^\alpha\,dx
=
\frac{\Gamma(n+\alpha+1)\Gamma(n+\beta+1)}
 {(2n+\alpha+\beta+1)\, n!\,\Gamma(n+\alpha+\beta+1)} \label{jacobi-norms}
.\end{multline}
For a function 
 $f(x)$ consider the numbers
 ('the Fourier coefficients') defined by the formula
$$
c_n(f):=\int_0^1 f(x)\,\cP_n^{\alpha,\beta}(x)\, x^{\alpha}(1-x)^{\beta}\,dx
.$$
Then the function
 $f(x)$  can be restored by the formula 
\begin{equation}
\label{inversion-j}
f(x)=\sum_{n=0}^\infty \frac { {c_n(f)}}{\gamma_n} \cP_n^{\alpha,\beta}(x)
.\end{equation}
Moreover, the following 'Plancherel formula' holds:
$$
\int_0^1 f(x)\,\ov{g(x)}\,\,x^\beta(1-x)^\alpha\,dx
=\sum_{n=0}^\infty \frac 1{\gamma_n} c_n(f)\,\ov{c_n(g)}
.$$

The expansion in the Jacobi polynomials has 
a continuous  analog, in which series are replaced by integrals.

\smallskip

{\bf\punct The index hypergeometric transform.}
Let $b$, $c>0$. For a function $f$
 on the half-line $[0,\infty)$ we define 
 a function of the variable
  $s\ge 0$ by
\begin{multline}
J_{b,c}f(s)=[\wh f]_{b,c}(s)
=\\
=\frac 1{\Gamma(b+c)}
\int_0^\infty f(x)  \F\Bigl[
  \begin{matrix} b+is,\, b-is\\b+c\end{matrix}
  ; -x\Bigr]\, x^{b+c-1} (1+x)^{b-c}\,dx
\label{def}
.\end{multline}

\begin{theorem}
\label{th:unitary}
 {\rm a)}
The operator  $J_{b,c}$ is a unitary operator  
$$
J_{b,c}:\,L^2\Bigl([0,\infty), \,  x^{b+c-1} (1+x)^{b-c}\,dx\Bigr)
\to L^2\Bigl([0,\infty),\densitybc ds\Bigr)
.$$

In other words, the following Plancherel formula holds
$$
\int_0^\infty f_1(x)\,\ov{f_2(x)}\densityx\,dx
=
\int_0^\infty [\wh f]_{b,c}(s) \ov{[\wh f]}_{b,c}(s)\densitybc\,ds
.$$

{\rm b)} The inverse operator is given by the formula 
\begin{equation}
f(x)=\frac 1{\Gamma(a+b)}
\int_0^\infty [\wh f]_{b,c}(s)
\F\Bigl[
  \begin{matrix} b+is,\, b-is\\b+c\end{matrix}
  ; -x\Bigr] \densitybc\, ds
\label{inversion}
.\end{equation}
\end{theorem}

Notice that  the statement
  b) follows from   a),
  because for a unitary operator
$U$  we have  $U^{-1}=U^*$.

As in the case of the Fourier transform,
we have  a question about a precise definition.
For instance, we can say that
the integral transform $J_{b,c}$ is defined 
on functions with compact support,
next we extend it by continuity to a unitary operator
defined in the space 
$L^2$.

\newpage

{\bf\punct Holomorphic extension to a strip.}

 \begin{lemma}
\label{l:holomorpic}
Let $f$ be integrable  on $\R_+$ and
$$f(x)=o(x^{-\alpha-\epsilon}), \qquad\qquad x\to+\infty,
$$
where $\epsilon>0$. Then $[\widehat f(s)]_{b,c}$
is holomorphic in the strip
$$|\Im s|<\alpha-b$$
and satisfy the condition
  $\wh f(-s)=\wh f(s)$ in the strip.
\end{lemma}

{\sc Proof.}
This follows from the asymptotics for the hypergeometric function,
(see \cite{HTF}, V. 1, (2.3.2.9)) as $x\to+\infty$:
$$
\FBC= \lambda_1 x^{-b+is} + \lambda_2 x^{-b-is}
+O(x^{-b+is-1}) +O(x^{-b-is-1})\, $$
where $2is\not\in \Z$, and $\lambda_1$, $\lambda_2$ are certain constansts 
(for  $2is\in \Z$ there arises  an additional factor $\ln x$ at
the leading term).
\hfill $\square$

\smallskip

{\bf\punct The operator calculus.}
 Denote by  $D$ the hypergeometric differential operator
\begin{equation}
D:=-x(x+1)\frac {d^2}{dx^2}-
\bigl[(c+b)+(2b+1)x\bigr] \frac {d}{dx}
+b^2
\label{D}
\end{equation}
(in comparision with the common notation
we  replaced  $x$  by $-x$).
The hypergeometric functions in
(\ref{def}) are  (generalized) eigenfunctions 
of the operator $D$:
 \begin{equation}
D  \F\Bigl[
  \begin{matrix} b+is,\, b-is\\b+c\end{matrix}
  ; -x\Bigr]=- s^2
\F\Bigl[
  \begin{matrix} b+is,\, b-is\\b+c\end{matrix}
  ; -x\Bigr]
\label{eigenfunctions}
.\end{equation}
It is easy to see that
  $D$ is formally self-adjoint in the following sense 
$$
\int_0^\infty Df_1(x)\cdot f_2(x)\,\densityx dx
=
\int_0^\infty f_1(x)\cdot D f_2(x)\,\densityx dx
$$
where $f_1$, $f_2$  are smooth compactly supported
functions on $(0,\infty)$
(in fact $D$ is essentially self-adjoint, see below)
and Theorem \ref{th:unitary}
is a theorem about expansion of  
 $D$ in eigenfunctions.

\begin{theorem}
\label{th:ss}
 Let $f$, $Df\in L^2$, then 
 \begin{equation}
[\wh{D f}(s)]_{b,c}=-s^2 \wh f(s)
.\end{equation}
\end{theorem}

{\sc Proof.} This is a rephrasing of formula
(\ref{eigenfunctions}).

 \begin{theorem}
\label{th:operation}
Let a function    $f$ be continuous on $\R^+$ 
and satisfies the condition
\begin{equation}
f(x)=o(x^{-b-1-\epsilon});\qquad x\to +\infty
\label{decreasing}
.\end{equation}
Then
\begin{equation}
[\widehat{xf(x)}]_{b,c}= P[\widehat{f(x)}]_{b,c}
,\end{equation}
where the difference operator
 $Pg$  is given by 
\begin{multline}
Pg(s)=\frac{(b-is)(c-is)}  {(-2is)(1-2is)}
    (g(s+i)-g(s))+
\frac{(b+is)(c+is)}  {(2is)(1+2is)}
    (g(s-i)-g(s))
.\end{multline}
\end{theorem}

{\sc  Remark.} Emphasize amusing characteristics
of this theorem.

\smallskip

1. The operator   $P$ is a difference operator,
but a shift 
  $s\mapsto s+i$ is in the imaginary direction,
  and integration is along the real axis.

\smallskip

2. The transformation $J^{-1}_{b,c}$  send the operator
 $P$ to the operator of multiplication
 by a function 
 $x$,  i.e., our operator   $J^{-1}_{b,c}$
 determines a spectral decomposition of the difference
 operator   $P$. For several examples of spectral decompositions
 of difference operators in imaginary direction, see 
\cite{Ner-difference}, \cite{Gro}.

\smallskip

3. The operator $P$ is similar to difference operators,
related to Wilson, continuous Hahn, continuous dual Hahn,
 Meixner--Pollachek orthogonal polynomials,
see \cite{AAR}, formulas
 (6.10.6), (6.10.9), (6.10.12) and Problem  6.37.c (see also \cite{KLS}).
 The rational coefficients of the operator
 $P$ are 'catenated' with the $\Gamma$-factors in the formula
(\ref{inversion}).

\smallskip

{\sc Proof.} This is reduced to a verification
of the identity 
 $$
P
\F\Bigl[
  \begin{matrix} b+is,\, b-is\\b+c\end{matrix}
  ; -x\Bigr]= x
\F\Bigl[
  \begin{matrix} b+is,\, b-is\\b+c\end{matrix}
  ; -x\Bigr]
.$$

\begin{theorem}
\label{th:operation-2}
 Let $f$ and $f'$ be continuous and decrease as
{\rm (\ref{decreasing})}. Then
\begin{equation}
[\widehat{x(x+1)\frac{d}{dx}f}]_{b,c}=H[\widehat f]_{b,c}
,\end{equation}
where the difference operator
 $H$ is given by 
\begin{multline}
Hg(s)
=
\frac{(b-is)(b+1-is)(c-is)}  {(-2is)(1-2is)}
         (g(s+i)-g(s))+\\
+
\frac{(b+is)(b+1+is)(c+is)}  {(+2is)(1+2is)}
         (g(s-i)-g(s))-(b+c)g(s)
.\end{multline}
\end{theorem}

{\sc Proof.}
We evaluate
 $J_{b,c}$-image of the commutator 
$[x, D]$.

\smallskip

{\bf\punct Historical remarks.}
The transformation
 $J_{1/2,1/2}$ was introduced by F.~G.~Mehler 
 \cite{Meh} in 1881г. He presented the inversion formula
 without proof  (it has to be said that the formula
 is not obvious at all). A proof was published
 by V.~A.~Fock
\cite{Fock} in 1943г. As a result,
the transformation
 $J_{1/2,1/2}$ is called the Mehler--Fock transform. 
The general transformation $J_{b,c}$ was introduced by H.Weyl
 in
 1910 the work \cite{Wey} on the spectral theory
 of of differntial operators. It seems
 that this result have not met the eye.
 Again this transform had appeared in the book
 of Titchmarsh   \cite{Tit} in 1946г.
 In 1949 the transform was rediscovered by
 M.~N.~Olevsky \cite{Ole}, 
apparently this was related to his works on
multi-dimensional Lobachevsky spaces. 

The most common terms for
 $J_{b,c}$ are the {\it Olevsky transform} and the {\it Jacobi transform} 
 (introduced by T.~Koornwinder).

\section{Applications to special functions}

\COUNTERS

Our first aim is to evaluate index hypergeometric transforms
of some functions.  We do this in Subsection
 2.2 by the Mellin transform.
 Next, in  2.3-2.4 we demonstrate effectivity
of the index transform as a tool of theory of special functions.

\smallskip

{\bf\punct The Mellin transform.}
Let $f(x)$ be a function defined on the half-line
 $x>0$. Its Mellin transform is defined by the
 formula 
$$
F(s)=\cM f(s):=\int_0^\infty f(x)\,x^{s}\frac {dx}x
.$$
The domain of absolute convergence of this integral is a
certain vertical strip of the form
$u<\Re s<v$, the function  $F(s)$ 
is holomorphic in this strip,
the boundaries of the strip can belong or do not belong
the domain of convergence; a strip can be degenerated
to a vertical line of the form  $\Re s=u$.
Certainly, it can be empty.
 
The Mellin transform is a unitary operator
from
 $L^2(\R,dx/x)$ to $L^2$ on vertical line  
 $\Re s=1/2$. In particular, 
$$
\int_0^\infty f(x)\ov{g(x)}\,dx/x=
\frac 1{2\pi}\int_{-\infty}^{+\infty}
F(1/2+is)\ov{G(1/2+is)}\,ds
.$$
Recall the theorem about convolution.
If the domains of definition of $F(s)=\cM f(s)$ and $G(s)=\cM g(s)$
have an intersection (a strip or a line),
then $\cM$ send the multiplicative convolution
$$
f*g(x):=\int_0^\infty f(y)g(x/y)\,dy/y
$$
to
 $F(s)G(s)$  (on the common domain of definition).

\smallskip

Certainly, the Mellin transform is reduced to the Fourier
transform by a substitution $x=e^y$.
From the point of view of abstract theory there is 
no difference between the Mellin transform and the Fourier
transform. But their role in  theory of special functions
is different.

\smallskip

{\bf\punct A game in the Mellin transform. A short table
of index transforms.}  Since we will meet long products
of $\Gamma$-functions, we will use the following notation

$$
\Gamma
\begin{bmatrix}a_1,&\dots,& a_k\\
               b_1,&\dots,& b_l  \end{bmatrix}
:=
\frac{\Gamma(a_1)\dots\Gamma(a_k)}
     {\Gamma(b_1)\dots\Gamma(b_l)}
.$$

Consider arbitrary convergent Barnes-type integral
$$
 \frac 1{2\pi i} \int_{-i\infty}^{+i\infty}
\Gamma
\begin{bmatrix}
a_1+s,\dots,a_k+s,b_1-s,\dots,b_l-s\\
c_1-s,\dots,c_m-s,d_1+s,\dots,d_n+s
\end{bmatrix}
x^{s}\,ds
.$$
It can be represented as a linear combination of hypergeometric
functions
 $_pF_q$  with $\Gamma$-coefficients. The idea is explained in the
 book \cite{AAR}, Section 2.4.
A calculation requires a watching of some asymptotics,
but it can be done once and forever in 'general case'.
 The final rules can be found in \cite{Sla}, \cite{PBM3}.

On the other hand there are unexpectedly many cases when
the integral admits a simpler expression than 
the result of the general algorithm,
see the tables of Prudnikov,  Brychkov, Marichev,
v.3, \cite{PBM3}, Chapter 8 (I do not know rational
explanations of this phenomenon).

Now we evaluate two auxiliary integrals.

\begin{lemma}
\label{l:mellin-1}
\begin{multline}
\int_0^\infty
\frac{x^{\alpha-1}}{(x+z)^\rho} \F(p,q;r;-x) \,dx=
\\=
\frac{z^{\alpha-\rho}} {2\pi i}
\Gamma\begin{bmatrix} r\\ p,q,\rho
\end{bmatrix}
\int_{-i\infty}^{i\infty}\Gamma
\begin{bmatrix} t+\alpha,\rho-t-\alpha,p+t,q+t,-t\\
r+t\end{bmatrix}
z^t   dt
\label{int1}
.\end{multline}

\begin{multline}
\int_0^\infty
x^{\alpha-1}
\F\left[\begin{matrix} p,q\\r \end{matrix}; -\omega x\right]
\F\left[\begin{matrix} u,v\\w \end{matrix}; -\tilde\omega x\right]
dx=\\=
\frac {\omega^{-\alpha}}{2\pi i}  \Gamma
 \begin{bmatrix} r,w\\ u,v,p,q\end{bmatrix}
\int\limits_{-i\infty}^{i\infty} \Gamma\begin{bmatrix}
\alpha+t,u+t,v+t,p-\alpha-t,q-\alpha-t,-t\\
 r-\alpha-t,w+t\end{bmatrix}
   \left(\frac{\omega}{\tilde\omega}\right)^{-t}dt
\label{int2}
.\end{multline}
\end{lemma}

{\sc Proof.} To be definite, we evaluate
the first integral. The Mellin transform of the function
$f(x):= x^{\alpha-1}/(x+z)^\rho$ is
$B(s+\alpha,\rho-s-\alpha)z^{s+\alpha-\rho}$.
The Mellin transform of
 $g(x):=\F(p,q;r;-x)$ is evaluated in \cite{AAR},
 Section 2.4,
and is a product of
$\Gamma$-functions. Our integral is a convolution
of $xf(1/x)$ and $g(x)$. Next, we observe that the Mellin transform
of the function  $xf(1/x)$ is $F(1-s)$.
It remains to apply the theorem about convolution.
 \hfill$\square$

\smallskip

Thus, our calculations was reduced to an rearrangement of 
 $\Gamma$-functions.

\smallskip

In the right hand side there are Barnes integrals,
which can be represented as linear combinations of functions
$\FF$
and   $\FFF$ respectively. We will not write them,
instead of this we notice that for some values of the
parameters $\Gamma$-factors in the right hand sides
can cancel.

\begin{lemma}
\label{l:mellin-2}
The transform $J_{b,c}$ send
\begin{align}
(1+x)^{-a-c}\qquad\ &\RA \qquad
\frac{\Gamma(c+is)\Gamma(c-is)}
     {\Gamma(c+a)\Gamma(c+b)}
\label{table1}
;
\\
 \frac{(1+x)^{b-a}}
       {(x+z)^{c+b}}
\qquad &\RA \qquad
     \Gamma \begin{bmatrix} c+is,c-is\\
     c+a,c+b \end{bmatrix}
\F\left[\begin{matrix} c+is,c-is\\c+a \end{matrix} ; 1-z\right]
\label{table2}
;\\
 x^{-u-a}   \qquad &\RA \qquad
\frac{\Gamma(-u+b)}{\Gamma(a+u)}
\cdot
\frac{\Gamma(u+is) \Gamma(u-is)}
      {\Gamma(b+is) \Gamma(b-is)}
\label{table3}
\end{align}
 \begin{multline}\F\left[\begin{matrix} p+b,q+b\\ a+b\end{matrix}
   ;- \dfrac xy \right](1+x)^{b-a}\RA
\\
y^{b-q}
\Gamma\begin{bmatrix}a+b\\p+q,p+b,q+b\end{bmatrix}
\cdot \Gamma\begin{bmatrix}p+is,p-is,q+is,q-is\\
a+is,a-is\end{bmatrix}
 \F \left[\begin{matrix} p+is,p-is\\p+q  \end{matrix}; 1-y \right]
\label{table4}
;\end{multline}
\begin{multline}
 \F\left[\begin{matrix} p+b,q+b\\ a+b\end{matrix}
   ;-x \right](1+x)^{b-a}
\RA
\\
\RA
 \Gamma\begin{bmatrix}a+b\\p+q,p+b,q+b\end{bmatrix}
\cdot \Gamma\begin{bmatrix}p+is,p-is,q+is,q-is\\
a+is,a-is\end{bmatrix}
\label{table5}
;\end{multline}
\begin{multline}
\F\left[\begin{matrix} a+c,a+d\\ a+b+c+d\end{matrix};-x\right]
\RA\\ \RA
\frac
{\Gamma(a+b+c+d)\cdot\Gamma(c+is)\Gamma(c-is)
\Gamma(d+is)\Gamma(d-is)}
{\Gamma(a+c)\Gamma(a+d)\Gamma(b+c)\Gamma(b+d)\Gamma(c+d)}
\label{table6}
.\end{multline}
\end{lemma}

{\sc Proof.} We look to the right hand side of 
(\ref{int2}). If $\alpha=r$, then two  
$\Gamma$-factors cancel. The rest is the integral
representation of 
$\F$. This gives the second formula. Substituting  $z=1$
we get the first formula.

Next, if
 $z=1$, $r=p+q+\rho$, then we get one of Barnes integrals
 in the right hand side 
(see Theorem  2.4.3 of \cite{AAR}). This gives us
(\ref{table3}).

Further, we watch possible simplifications in the right hand side of
(\ref{int2}). After the substitution   $\alpha=w=r$
we get a cancelation of four $\Gamma$-factors.
This gives formula
 (\ref{table4}). Substiting $y=1$
to (\ref{table4}) we get (\ref{table5}).

To verify (\ref{table6}),
we observe that two   $\Gamma$-factors in the right
hand side of 
 (\ref{int2}) cancel, and we again apply Theorem  2.4.3
 of \cite{AAR}.

\smallskip

{\bf\punct Game in the Plancherel formula.}
Thus we wrote a short table with 6 row for
the transform  $J_{b,c}$.
Applying the Plancherel formula
for
 $J_{b,c}$ we can get an amusing collection
 of integrals. We present several examples. 

\smallskip

a) {\it The De Branges--Wilson integral.}
Applying the Plancherel formula for pair of functions
$(1+x)^{-a-c}$ and $(1+x)^{-a-d}$,
we after a trivial calculation we get 
the De Branges--Wilson integral
  (see Section 3.6 of \cite{AAR}), 
  this derivation is due to
  Koornwinder,\cite{Koo3}. Recall that its is given by
\begin{equation}  
  \frac 1 {2\pi}
\int_{-\infty}^\infty
\left|\frac{\prod_{k=1}^4
            \Gamma(a_k+is)}
     {\Gamma(2is)}\right|^2 ds=
     \frac{\prod_{1\le k< l\le 4}\Gamma(a_k+a_l)}
     {\Gamma(a_1+a_2+a_3+a_4)}
     \label{eq:de-branges}
     .
  \end{equation}

\smallskip

b) {\it Another beta-integral.} 
Applying the Plancherel formula to
   $x^{-u-a}$, $x^{-v-a}$,
we get the integral
$$ \frac 1 {2\pi}
\int_{-\infty}^\infty
\left|\frac{\prod_{k=1}^3
            \Gamma(a_k+is)}
     {\Gamma(2is)\Gamma(b+is)}\right|^2 ds=
\frac{\Gamma(b-a_1-a_2-a_3)
\prod_{1\le k< l\le 3}\Gamma(a_k+a_l)}
 {\prod_{k=1}^3 \Gamma(b-a_k)}
.$$

c) {\it  An integral representation for ${}_3F_2(1)$.}
  The pair of functions
$$
(1+x)^{-a-e}\qquad\text{and} \qquad
\F\left[\begin{matrix} a+c, a+d\\a+b+c+d\end{matrix};-x\right]
$$
gives the following integral representation of
$\FF(1)$,        
\begin{multline}
\frac{1}{\pi}\int_0^\infty
\left|\frac{\Gamma(a+is)\Gamma(b+is)\Gamma(c+is)\Gamma(d+is)
\Gamma(e+is)}
{\Gamma(2is)}\right|^2ds=\\=
\frac{\Gamma(a+b)\Gamma(a+c)\Gamma(a+d)\Gamma(a+e)
\Gamma(b+c)\Gamma(b+d)\Gamma(b+e)\Gamma(c+d)\Gamma(c+e)}
{\Gamma(a+b+c+d)\Gamma(a+b+c+e)}
\times\\ \times
 \FF\left[\begin{matrix} a+c, b+c, a+b\\ a+b+c+d, a+b+c+e\end{matrix}
;1\right]
.\end{multline}
The left hand side is symmetric with respect to the parameters,
therefore the right hand side also is symmetric.
This symmetry is the Kummer identity
(see \cite{AAR}, Corollary 3.3.5).

\smallskip

d) {\it Adding a $\Gamma$-factor to the numerator.}
Applying the Plancherel formula to the pair 
 $$
  \F\left[\begin{matrix} a+c, a+d\\a+b+c+d\end{matrix};-x\right]
\qquad\text{ and} \qquad
   \F\left[\begin{matrix} a+e, a+f\\a+b+e+f\end{matrix};-x\right]
,$$
we get the identity
\begin{multline}
\frac 1\pi
\int_0^\infty
  \left|\frac{\Gamma(a+is)\Gamma(b+is)\Gamma(c+is)\Gamma(d+is)
\Gamma(e+is)\Gamma(f+is)}
{\Gamma(2is)}\right|^2ds=\\ =
\frac 1{2\pi i}
\Gamma(a+c)\Gamma(a+d)\Gamma(c+d)
\Gamma(b+e)\Gamma(b+f)\Gamma(e+f)\times
\\ \times
\frac1{2\pi i}\int_{-i\infty}^{+\i\infty}
\Gamma\begin{bmatrix} a+b+s, a+e+s, a+f+s, d-a-s,c-a-s,-s
      \\
       c+d-s,a+b+e+f+s
       \end{bmatrix}
ds
.\end{multline}
The right hand side is a linear combination of three functions
 $_4F_3(1)$ with $\Gamma$-coefficients.
 By the way a Barnes integral can be regarded as a final answer.

\smallskip

e) {\it Adding a $\Gamma$-factor to the denominator.}   
Now we apply the Plancherel formula to the pair
$$
\F\left[\begin{matrix} p+b,q+b\\a+b\end{matrix};-x\right]
    (1+x)^{b-a}
\quad\text{and} \quad
 \F\left[\begin{matrix} u+b,v+b\\a+b\end{matrix};-x\right]
    (1+x)^{b-a}
.$$
We omit intermediate calculations and present the final result
 \begin{multline}
\frac 1\pi
\int_0^\infty
\left|\frac{\Gamma(b+is)\Gamma(p+is)\Gamma(q+is)\Gamma(u+is)
\Gamma(v+is)}
{\Gamma(2is)\Gamma(a+is)}\right|^2ds =\\=
\frac 1{2\pi i}
\Gamma\begin{bmatrix} u+v,p+q,p+b,q+b\\ a-v,u-v\end{bmatrix}
\times\\ \times  \int_{-i\infty}^{i\infty}
 \Gamma\begin{bmatrix} u+p+s,u+q+s,b+u+s,a-v+s,v-u-s,-s\\
u+a+s,u+b+p+q+s\end{bmatrix} 
     .
\label{eq:6-2}
\end{multline}

Two last identities are not as aesthetic
as previous. However, consider two special cases of the last integral.

\smallskip

f) {\it The Nassrallah--Rahman integral.}
In the last integral we set
$a=b+u+v+p+q$ (this leads to cancellation of $\Gamma$-factors)
and applying Theorem 2.4.3 of \cite{AAR}
 (after changing notation),
 we obtain the Nassrallah--Rahman integral
(its $q$-version is present in the book \cite{AAR},
Theorem 10.8.2)
$$
\frac1{2\pi}\int_{-\infty}^\infty
 \left| \frac{\prod_{j=1}^5 \Gamma(a_j+is)}
   {\Gamma(2is)\Gamma(\sum_1^5 a_j+is)}\right|^2ds=
2\frac{\prod_{1\le k<l \le 5}\ \Gamma(a_k+a_l)}
  {\prod_{k=1}^5\Gamma(a_1+a_2+a_3+a_4+a_5-a_k)}
.$$

g) {\it The Whipple identity 
and the symmetry of the Wilson polynomials%
\footnote{Recall that the Wilson polynomials $W_n(a,b,c,d;s^2)$
 are given by the formula
$$
W_n(a,b,c,d;s^2)=(a+b)_n (a+c)_n (a+d)_n
\,\,
\FFF\left[\begin{matrix}
-n,n+a+b+c+d-1, a+is,a-is\\
a+b,a+c,a+d\end{matrix}; 1\right].
$$
They are orthogonal with respect to the weight
(see (\ref{eq:de-branges})
 $$w(s) =
 \frac 1\pi\left|\frac{\Gamma(a+is)\Gamma(b+is)
 \Gamma(c+is)\Gamma(d+is)}
{\Gamma(2is)}\right|^2.$$}
 with respect to the parameters.}   
We represent the right hand side of  (\ref{eq:6-2})
in terms of $_4F_3$,
\begin{multline}
\Gamma\begin{bmatrix}
u+v,p+q,p+b,q+b\\ a-v,u-v
\end{bmatrix}\times\\ \times
\Biggl\{
\Gamma\begin{bmatrix} v-u,u+p,u+q,u+b,a-v\\
  u+a,u+b+p+q\end{bmatrix}
\FFF\left[\begin{matrix} u+p,u+q,u+b,a-v\\
1+u-v,u+a,u+b+p+q\end{matrix};1\right]
+\\+
\Gamma\begin{bmatrix} u-v,p+v,q+v,b+v,a-u\\
       v+a,v+b+p+q\end{bmatrix}
 \FFF\left[\begin{matrix}
     p+v,q+v,b+v,a-u\\
     1-v+u,v+a,v+b+p+q
      \end{matrix};1\right]
      \Biggr\}
.
\label{eq:6-add}
\end{multline}
The left hand side of
 (\ref{eq:6-2}) is symmetric with respect to the parameters
$b$, $p$, $q$, $u$, $v$, but the right hand side in the form 
(\ref{eq:6-add}) does not looks as symmetric.
This gives a symmetry relation for
 $_4F_3$ in the form 
"a linear combination of four summands equals  0".
This is the  "nonterminating Whipple identity". 
Its "terminating version"  (Theorem 3.3.3 of \cite{AAR})
can obtained by the substitution 
 $a=v-m$ with integer  $m$, then two summands disappear
 due factors  $\Gamma(-m)$ in the denominators 
(this identity is 
the symmetry of the Wilson polynomials with respect to the parameters).

\smallskip

h) {\it One again extension of the De Branges--Wilson integral.}
Applying the Plancherel formula  to
$ (1+x)^{b-a}(1+x+y)^{-b-c}$  and
 $ (1+x)^{b-a}(1+x+y)^{-b-d}$, we get
\begin{multline}
 \frac 1\pi \int_0^\infty
\left| \frac{\Gamma(a+is)\Gamma(b+is)\Gamma(c+is)\Gamma(d+is)}
            {\Gamma(2is)}\right|^2
\times\\ \times
\F\left[\begin{matrix} c-is,c+is\\a+c\end{matrix}; -y\right]
 \F\left[\begin{matrix} d-is,d+is\\a+d\end{matrix}; -y\right] \,ds
=\end{multline}
$$
 =\frac{\pi\Gamma(a+b)\Gamma(a+c)\Gamma(a+d)\Gamma(b+c)
\Gamma(b+d)\Gamma(c+d)}{\Gamma(a+b+c+d)}
\F\left[\begin{matrix} 2b+c+d,c+d\\a+b+c+d\end{matrix};-y
\right]
. $$

{\bf \punct A derivation of the orthogonality relations
for Wilson polynomials.}%
\footnote{Another simple derivation is contained in
 \cite{Ner-wilson}.}
 Now we derive
  $J_{b,c}^{-1}$-image of the function 
$$|\Gamma(a+is)|^2 W_n(s^2), \quad\text{
where $W_n(s^2)=W_n(a,b,c,d;s^2)$ is a Wilson polynomial}.
$$
We must evaluate the integral
\begin{multline*}
\frac 1{\Gamma(b+c)}\int_0^\infty
\left| \frac{\Gamma(a+is)\Gamma(b+is)\Gamma(c+is)}
{\Gamma(2is)}\right|^2
\F\left[\begin{matrix}b+is,b-is\\b+c\end{matrix};-x\right]
W_n(s^2)\,ds
=\\=
\frac {(a+b)_n(a+c)_n(a+d)_n}{\Gamma(b+c)}
\int_0^\infty
\left| \frac{\Gamma(a+is)\Gamma(b+is)\Gamma(c+is)}
{\Gamma(2is)}\right|^2
\times\\ \times
\F\left[\begin{matrix}b+is,b-is\\b+c\end{matrix};-x\right]
\sum_{k=0}^n
\frac{(-n)_k(n+a+b+c+d-1)_k(a+is)_k(a-is)_k}
{k!(a+b)_k(a+c)_k(a+d)_k}\,ds
.\end{multline*}
We get a linear combination of known for us
(in virtue of the inversion formula and (\ref{table1}))
integrals of the type 
\begin{multline*}
\int_0^\infty
\left| \frac{\Gamma(a+k+is)\Gamma(b+is)\Gamma(c+is)}
{\Gamma(2is)}\right|^2
\F\left[\begin{matrix}b+is,b-is\\b+c\end{matrix};-x\right]
\,ds
=\\=\frac{\Gamma(a+b+k)\Gamma(a+c+k)}{(1+x)^{a+b}}
.\end{multline*}
As a result we get
$$
\frac{\Gamma(a+b)\Gamma(a+c)(a+b)_n}
{\Gamma(b+c)}(1+x)^{-a-b}
\F\left[\begin{matrix}-n,n+a+b+c+d-1\\a+d\end{matrix}
;\frac 1{1+x}\right]
.$$
This is a Jacobi polynomial of the variable
$1/(1+x)$.

\smallskip

Transposing
 $a$ c $d$, we evaluate the inverse index transform
 of 
$|\Gamma(d+is)|^2W_m(s^2)$.
Next we evaluate integral
 (see formula (3.8.3) of \cite{AAR})
$$
\frac 1{\Gamma(b+c)}\int_0^\infty
\left| \frac{\Gamma(a+is)\Gamma(b+is)\Gamma(c+is)\Gamma(d+is)}
{\Gamma(2is)}\right|^2
W_n(s^2)\,W_m(s^2)\,ds
$$
applying the Plancherel formula.
We get
\begin{multline*}
\const \int_0^\infty
\F\left[\begin{matrix}-n,n+a+b+c+d-1\\a+d\end{matrix}
;\frac 1{1+x}\right]
\times\\ \times
\F\left[\begin{matrix}-m,m+a+b+c+d-1\\a+d\end{matrix}
;\frac 1{1+x}\right] x^{b+c-1}(1+x)^{-a-d-c} dx
.\end{multline*}
Passing to the variable
 $y=1/(1+x)$, we get the integral, which express the orthogonality
 relations for the Jacobi polynomials.

At the first glance this proof of orthogonality can be invented
only if we know the final result.
Really this calculation gives us the following observation:

{\it The Wilson orthogonal system with  $a=d$ 
is the image of the Jacobi system under the index transform} (\cite{Koo3}). 

Notice that the index transform was discovered in 
1910 and became well-known upto 1950,
therefore it seems strange that the Wilson polynomials were discovered so late
(1980).

\section{Derivation of the inversion formula. Jump of resolvent}

\COUNTERS
 
Many ways of derivations are known, see \cite{Koo2}.
In particular, we can decompose the index transform
as a products of simpler integral transforms 
and apply inversion formulas for the factors.
However, the original way of Weyl based 
on spectral theory seems the most natural
up to now
(see, e.g., \cite{DS}, \S13.8, or \cite{Tit}).
We present a version of derivation
using a minimum of theory but requiring
superfluous calculations.
 In details, the spectral theory of differential
 operators is exposed in
Titchmarsh \cite{Tit}, Dunford, Schwartz \cite{DS}, Chapters
XII--XIII, and Naimark \cite{Nai}.

\smallskip

{\bf \punct Jump of the resolvent.}
Recall the Spectral Theorem.
Consider a finite or countable collection 
of measures   $\mu_1$, $\mu_2$,\dots
on $\R$, the Hilbert space  
$$V[\vec \mu]:=\oplus_j L^2(\R,\mu_j),$$
and an operator $Z_{\vec\mu}:V[\vec\mu]\to V[\vec\mu]$
given by the formula
$$
\bigl[Z_{\vec\mu}(f_1\oplus f_2\oplus\dots)\bigr](x)=
xf_1(x)\oplus xf_2(x)\oplus\dots
$$

\begin{theorem}
\label{th:abstract}
For any self-adjoint (generally, unbounded) operator in a Hilbert space
 $H$ there exists a collection of 
  $\mu_j$  and a unitary operator  $U:H\to V[\vec\mu]$
such that $A=U^{-1}Z_{\vec\mu}U$
\end{theorem}

For any Borel subset
 $M\subset \R$ consider the subspace 
  $W(M)\subset V[\vec\mu]$ of functions,
  which equals 0 outside 
the set $M$. Define the spectral subspace
$\Omega(M):=U^{-1}W(M)$.
Denote by $P[\Omega]$ the projection operator
to this subspace.

\begin{proposition}
\label{pr:}
 For any finite interval  $(a,b)\subset\R$,
$$P\bigl[(a,b)\bigr]=
\frac1{2\pi i}
\lim_{\delta\to+0}\lim_{\epsilon\to+0}
\int_{a+\delta}^{b-\delta}\bigl( (\lambda-i\epsilon-A)^{-1}
-(\lambda+i\epsilon-A)^{-1}\bigr)\,d\lambda
.
$$
\end{proposition}

The limit here is the limit in the strong operator topology,
 $T_n\to T$ if for any vector
 $v$ we have  $\|T_nv-Tv\|\to 0$.

A verification of this statement is straighforwarward
(and is a good exercise, in particular for
finite-dimensional spaces),
we can from  outset assume that our operator acts
in  $V[\vec\mu]$.

For any vector
 $v$,
$$
v=
\frac1{2\pi i}
\lim_{N\to\infty}\lim_{\epsilon\to+0}
\int_{-N}^{N}\bigl( (\lambda-i\epsilon-A)^{-1}
-(\lambda+i\epsilon-A)^{-1}\bigr) v\,d\lambda
.$$

An evaluation of the limit gives the spectral decomposition.
Now we will perform this for the hypergeometric  differential operator 
 $D$ defined above
 (\ref{D}).

\smallskip

{\bf\punct Solutions of the equation $(D-\lambda)f=0$.}
For each $\lambda$ this equation has two linear
independent solutions. We choose two bases
in the space of solutions
(both bases consist of Kummer series, see \cite{HTF}, Section 2.9).
The first basis consists of functions
\begin{align}
\phi(x,\lambda)&=\F[b+\sqrt\lambda,b-\sqrt\lambda;b+c;-x];\\
\psi(x, \lambda)&=(-x)^{1-b-c}\F[1+\sqrt\lambda-c,1-\sqrt\lambda-c;2-b-c;-x]
.\end{align}
The second basis $u_\pm(x)$  is given by formulas
\begin{equation}
u_\pm(x, \lambda)=(-x)^{-b\mp\sqrt\lambda}
\F[b\pm\sqrt\lambda,1\pm\sqrt\lambda-c;1\pm 2\sqrt\lambda;-x^{-1}]
.\end{equation}
We assume that the complex plane 
 $\lambda$ is cut along the negative semi-axis.

For the first pair of functions the behavior near zero
is easily observable, for the second pair we see the behavior near infinity.
Below we need a formula expressing 
 $\phi$ in terms of
$u_+$ and $u_-$:
$$
\phi(x, \lambda)=B_+(\lambda)u_+(x, \lambda)+
B_-(\lambda)u_-(x,\lambda)
,$$
where
\begin{equation}
B_\pm(\lambda)=\frac{\Gamma(b+c)\Gamma(\mp\sqrt \lambda)}
                     {\Gamma(b\mp \sqrt\lambda)\Gamma(c\mp\sqrt\lambda)}
\label{Bpm}
.\end{equation}

{\bf\punct Self-adjointness.} Let $b>0$, $c>0$.
We define the operator
$D$ on the space  $\cD(\R_+)$ of smooth compactly supported functions
on $(0,\infty)$. The operator  $D$ is formally symmetric with respect to the weight
  $x^{b+c-1}(1+x)^{b-c}\,dx$,
i.e.,
$$
\int_0^\infty (Df)(x)\ov{g(x)}x^{b+c-1}(1+x)^{b-c}\,dx=
  \int_0^\infty f(x)\ov{Dg(x)}x^{b+c-1}(1+x)^{b-c}\,dx
$$
where $f$, $g\in\cD(\R_+)$. Its adjoint operator $D^*$
is determined from the condition  $D^*g=h$ if $g$,
$h\in L^2(\R_+,x^{b+c-1}(1+x)^{b-c})$ and
$$
\int_0^\infty (Df)(x)\ov{g(x)}x^{b+c-1}(1+x)^{b-c}\,dx=
  \int_0^\infty f(x)\ov{h(x)}x^{b+c-1}(1+x)^{b-c}\,dx
$$
for all $f\in\cD(\R_+)$. As before, this operator is given by
formula 
 (\ref{D}), but its domain of definition have increased.

Recall that for any formally symmetric operator
$A$ the numbers $\dim\ker(A^*-\lambda)$ (deficiency indexes)
are constant on the half-planes 
$\Im \lambda>0$ and $\Im\lambda<0$. The operator  $A$ 
is essentially self-adjoint if the both numbers are 0.
Therefore we must verify existence/nonexistence of solutions of
the differential equations $Df=\lambda f$ with $\Im\lambda\ne 0$
such that $f$ is contained in
 $L^2$ with respect to our weight. To be definite consider
 the upper half-plane  $\Im\lambda>0$.

It is easy to see that for
 $b+c>2$ such solutions do not exist.
Indeed, $\psi$ is too large near 0,
and $u_-$ is too large near $\infty$.
Therefore an $L^2$-solution must coincide with $\phi$ and $u_+$
simultaneously. But these two solutions are different.

Therefore  $D$ is essentially self-adjoint.

\smallskip

{\sc Remark.}
If  $b+c<2$, then $\phi$ and $\psi$ are in  $L^2$ near  
0. Therefore $u_+\in L^2$ and the operator $D$ is not self-adjoint.
We extend the operator
 $D$ and define it on the space of functions
 smooth on the closed  half-line $[0,\infty)$
 and vanishing for large
 $x$.
Then the operator became self-adjoint.
Below we do not watch this case.
\hfill $\square$

\smallskip

{\bf\punct The resolvent.}

\begin{lemma}
\label{l:resolution}
The resolvent $(D-\lambda)^{-1}$ of the operator  $D$
is defined in the domain 
 $\C\setminus [-\infty,0)$ and is given by 
\begin{equation}
L(\lambda)f(x)=\int_0^\infty K(x,y;\lambda)
 y^{b+c-1}(1+y)^{b+c}\,dy
\label{resolvent}
,\end{equation}
where the 'Green function' $K$ is given by
\begin{equation}
K(x,y;\lambda)=\begin{cases}
2B_-(\lambda)^{-1}\lambda^{-1/2} \phi(x,\lambda) u_+(y,\lambda),\qquad
\text{if $x\le y$}\\
2B_-(\lambda)^{-1}\lambda^{-1/2}  \phi(y,\lambda) u_+(x,\lambda),\qquad
\text{if $x\ge y$}\\
\end{cases}
\label{K}
,\end{equation}
 $B_-(\lambda)$ is defined by {\rm(\ref{Bpm})}.
\end{lemma}

The jump of the resolvent appear on the semi-axis
 $\lambda\le 0$ due discontinuity of 
 $\sqrt\lambda$ on the cut.
 Evaluating the jump of resolvent we get
\begin{equation}
f(x)=\frac 1{2\pi} \int_{-\infty}^0 \frac{d\lambda}{2\sqrt{\lambda}}
 \frac{\phi(y,\lambda)}{B_+(\lambda)B_-(\lambda)}
 \int_0^\infty \phi(z,\lambda)f(z)z^{b+c-1}(1+z)^{b-c}\,dz
.
\label{eq:inversion-bis}
\end{equation}

{\sc Remark.} Formula is so simple, because 
 $\phi$ has no jump; jump of
 $u_+$ is proportional to $\phi$.
 \hfill $\square$

\smallskip 
 
 The last formula is the desired inversion formula.

\smallskip

{\sc Proof of Lemma.} First, we formally check the identity
 $(D-\lambda)L(\lambda)=1$.
 We must verify that the function
$K$ satisfies the equation 
\begin{equation}
(D_x-\lambda) K(x,y;\lambda) y^{b+c-1}(1+y)^{b+c}=\delta(x-y)
\label{K2}
.\end{equation}
Obviously, outside the diagonal
 $x=y$ the equality 
$(D_x-\lambda) K=0$ holds. The kernel $K$
is continuous, but the first derivative has a jump.
Therefore,
\begin{multline*}
(D_x-\lambda) K(x,y;\lambda)=x(x+1)
\Bigl\{\frac
{\partial K(x,y,\lambda)}{\partial x}\Bigr|_{y=x+0}-
\frac{\partial K(x,y,\lambda)}{\partial x}\Bigr|_{y=x-0}
\Bigr\}\delta(x-y)
=\\=
2B_-(\lambda)^{-1}\lambda^{-1/2}
\Bigl[\phi(y,\lambda)'u_+(y,\lambda)-
 \phi(y,\lambda) u_+(y,\lambda)'\Bigr]\delta(x-y)
.\end{multline*}
In square brackets we have Wronski determinant of two solutions of
 $(D-\lambda)f$.
 Upto a  constant factor Wronskian is detemined by a differential equation,
 in our case it equals
$\const\cdot y^{-b-c}(1+y)^{c-b-1}$.
To evalute the constant factor we watch asymptotics 
of the Wronskian as 
$y\to\infty$.

In fact this calculation is sufficient to a proof.
But a priory boundedness of $L(\lambda)$ in $L^2$
is not evident. We overcome this difficulty
in the following way.

Since 
 $D$ is essentially self-adjoint, for $\lambda\notin\R$
the operator $(D-\lambda)^{-1}$ is unbounded. 
In virtue of L.~Schwartz's Kernel Theorem
(see, e.g.,
 (см. \cite{Her}) $(D-\lambda)^{-1}$
 is an integral operator, its kernel
 $K(x,y;\lambda)$ is a distribution of two variables.
  It satisfies  equation
 (\ref{K2}) and the symmetry condition
$K(y,x;\lambda)=\ov{K(x,y,\ov\lambda)}$.
Therefore, outside the diagonal 
 $x=y$ the distribution satisfies the system of equations
$$(D_x-\lambda)K=0,\qquad (D_y-\lambda)K=0.$$
It can be readily checked that our kernel
 $K$ is a unique admissible candidate,
 all other solutions of the system increase too rapidly.
 \hfill $\square$

\smallskip

{\bf\punct The  Romanovski Polynomials.}
Now let
$b<0$, $b+c>0$. Let $m=0,1,\dots, [-b]$.
 Consider the polynomials $p_m$ given by
\begin{multline*}
p_m(x):=
\F\left[\begin{matrix}-m,2b+m\\b+c\end{matrix};-x
     \right]=\\=
\frac{x^m\Gamma(b+c)\Gamma(-m-b)}
                          {\Gamma(2b+m)\Gamma(c+b+m)}
  \F\left[\begin{matrix}-m,1-m-b-c\\
2-b-c\end{matrix};-\frac 1x\right]
.
\end{multline*}

\begin{theorem}
\label{th:romanovsky}
{\rm a)} The polynomials $p_m$ are contained in
$L^2(\R_+,x^{b+c-1}(1+x)^{b-c})$.

\smallskip

{\rm b)} $Dp_m=(b+m)^2 p_m$.

\smallskip

{\rm c)} The polynomials $p_m$ are pairwise orthogonal.
\end{theorem}

The statements
 a), b) are evident,  с) follows from a) and b).
 \hfill $\square$

\smallskip                                      

Thus we get a finite system of orthogonal polynomials.
We can not enlarge it because the monomials
 $x^N$ with larger powers are not in 
 $L^2$.

\smallskip

{\sc Remark.}
A lot of such finite systems of orthogonal polynomials
is known.
Romanovski \cite{Rom} introduced another two systems:
the polynomials on the line orthogonal with respect to the weight
$$\frac{dx}{(1+ix)^\alpha(1-ix)^{\ov \alpha}},\qquad \text{where $\alpha\in\C$,}$$
(they also are analytic continuations
of the Jacobi polynomials)
and polynomials on the half-line $x\ge 0$
with respect to the weight
$$x^{-\beta} \exp(-1/x)\,dx$$
 (this is an analytic continuation
 of the Laguerre polynomials of  $1/x$).
 More complicated finite systems of orthogonal polynomials were
 enumerated by P.~Lesky
 \cite{Les1}, \cite{Les2}, some additions are in
 \cite{Ner-wilson}).  
 \hfill $\square$
 
 \smallskip

Our considerations explain this phenomenon.
For $b<0$, $b+c>0$  our operator  $D$
has finite number of discrete eigenvalues corresponding
to the Romanovski polynomials, which are 
added to the continuous spectrum.
 
It is necessary to modify our calculation.
The resolvent   (see formula (\ref{K})) now has
a finite number of poles at points
$\lambda=(b+m)^2$, they arise from the poles of 
 $B_-(\lambda)^{-1}$. To write the jump of the resolvent
we must additionally evaluate residues at these poles.
In the  inversion formula  (\ref{eq:inversion-bis}) we get additional terms
$$
\dots
+\sum_m \frac{\la f, p_m\ra_{L^2}} {\la p_m, p_m\ra_{L^2}} p_m(x).
$$

The expression (\ref{def}) for $J_{b,c}$ does not change.
But the function   
$J_{b,c}f(s)$ now is defined on the following
subset in $\C$:
the half-line $s\ge 0$ and a finite set  $s=i(b+m)$
on the imaginary axis (these points correspond to Romanovski polynomials).

\smallskip
 
{\sc Remark.} Such  orthogonal systems arise in
non-commutative harmonic analysis and correspond to 
discrete part of  spectra
(for instance, for  $L^2$ on pseudo-Riemannian 
symmetric spaces of rank 1,
 see also \cite{Ner-hahn}
 about tensor products of unitary representations of the group
 $\SL(2,\R)$.
Apparently (nobody verified this)
the discrete Flensted-Jensen series
   \cite{Fle} are controlled by some 
   multivariable orthogonal systems of Romanovski type.

\section{Applications to harmonic analysis}

\COUNTERS

{\bf\punct Pseudounitary groups of rank 1.} 
Let  $\K$ be $\R$, $\C$
or the quaternion algebra
 $\H$. 
 The case $\K=\R$ is sufficiently interesting.
We present several simple facts without proofs,
the reader can believe or verify.

Denote by $r$ the dimension of  $\K$.
 Let  $\K^n$ be the 
$n$-dimensional space over $\K$  with the standard
inner product,
$$\langle z,u\rangle=\sum  z_j \overline u_j
.$$
By $\U(1,n;\K)$ we denote the  {\it pseudounitary group} over $\K$,
i.e., the group of
  $(1+n)\times(1+n)$-matrices
$\bigl(\begin{smallmatrix}a&b\\c&d\end{smallmatrix}\bigr)$ over $\K$
satisfying the condition
$$
\begin{pmatrix}a&b\\c&d\end{pmatrix}
\begin{pmatrix}-1&0\\0&1\end{pmatrix}
\begin{pmatrix}a&b\\c&d\end{pmatrix}^*=
\begin{pmatrix}-1&0\\0&1\end{pmatrix}
.$$
The standard notations for the groups 
 $\U(1,n;\K)$
for $\K=\R,\C,\H$ are respectively:
$\OO(1,n)$, $\U(1,n)$, $\Sp(1,n)$.

\smallskip

{\bf\punct Homogeneous hyperbolic spaces.}
Denote by $\B_n(\K)$ the open unit ball 
$\langle z,z\rangle <1$ in $\K^n$. By $S^{rn-1}$ we denote the sphere 
$\langle z,z\rangle=1$. The group  $\U(1,n;\K)$ acts on $\B_n(\K)$
by linear-fractional transformations 
\begin{equation}
z\mapsto z^{[g]}:=(a+zc)^{-1}(b+zd)
.\end{equation}
The stabilizer 
 $K$  of the point  $0\in \B_n(\K)$
 consists of matrices of the form
\begin{equation}
\begin{pmatrix}a&0\\0&d\end{pmatrix}
                     \qquad\qquad |a|=1,\quad d\in\U(n;\K)
.\end{equation}
Therefore  $\B_n(\K)$  is the homogeneous space 
$$\B_n(\K)=\U(1,n;\K)\bigl/(\U(1;\K)\times\U(n;\K))
.$$

\smallskip

{\sc Remark.}
If   $\K=\R$, then our ball is the
 $n$-dimensional Lobachevsky space in the Beltrami--Klein  model.
 Recall that in this case straight lines in the Lobachevsky sense
 are segments (chords), the sphere 
$S^{n-1}$ is the absolute in the Lobachevsky sense.
The group $\OO(1,n)$ is the group of motions
of the Lobachevsky space. 

For $\K=\C$ and $\H$ we get complex and quaternionic
hyperbolic spaces. 
\hfill$\square$

\smallskip

The Jacobian of the transformation
 (1.5) is
$$J(g;z)=|a+zc|^{-r(1+n)}
.$$

Note the following simple formula
$$1-\langle z^{[g]},u^{[g]} \rangle
=(a+zc)^{-1}(1-\langle z,u\rangle)\overline{ (a+uc)}\,^{-1}
.$$
This implies than the
 $\U(1,n;\K)$-invariant measure on  $\B_n(\K)$
has the form
$$dm(z)=(1-\langle z,z\rangle)^{-(n+1)r/2} dz,$$
where $dz$ denotes the Lebesgue measure on    $\B_n(\K)$.

The group
$\U(1,n;\K)$
acts in  $L^2(\B_n(\K),  dm(z))$ by changes of variable 
\begin{equation}
\rho(g)f(z)= f\bigl((a+zc)^{-1}(b+zd) \bigr)
.\end{equation}
Evidently these operators are unitary.
In other words we get an infinite dimensional unitary representation
of the group
 $\U(1,n,\K)$.

Our next problem is {\it to decompose this representation into irreducible
representations}.

\smallskip

{\bf\punct The spherical principal series.}
Let $s\in\R$.
A representation
 $T_s$ {\it of spherical principal series}
 of the group  $\U(1,n;\K)$
is realized in  $L^2(S^{rn-1})$ and
is given by the formula
\begin{equation}
T_s
\begin{pmatrix}a&b\\c&d\end{pmatrix}
f(h)=f\bigl( (a+hc)^{-1}(b+hd)\bigr)|a+hc|^{-(n+1)r/2+1+is}
,\end{equation}
where $h\in S^{rn-1}$.
A straightforward calculation shows that these representations are unitary
for $s\in\R$.

\smallskip

{\sc Remark.}
All representations $T_s$ are irreducible,
representations
 $T_s$  and $T_{-s}$ are equivalent (this is not completely
 obvious, see  \cite{Vil}).

\smallskip

{\sc Remark.}
The term 'series' is used because these groups
have different types of unitary representations.
The term 'spherical' means that any representation $T_s$
contains a (unique) $K$-invariant vector.
In our model this vector is the function 
 $f=1$.

\smallskip

{\bf \punct  An intertwining operator.}
Consider the space of functions 
 $\phi(h,s)$ on the semi-cylinder  $S^{rn-1}\times \R_+$
(a precise description of this space is given below), let 
$\U(1,n;\K)$ act in this space by the formula 
$$
\tau\begin{pmatrix}a&b\\c&d\end{pmatrix}\phi(h,s)=
\phi\bigl( (a+hc)^{-1}(b+hd),s\bigr)\,|a+hc|^{-(n+1)r/2+1+is}
.$$
For a fixed 
   $s$ we get the representation  
$T_s(g)$ in functions depending on $h$. Thus we have some kind of a direct sum
of all representations  $T_s$ with respect to a continuous parameter  $s$
(thus is called a 'direct integral').

We define the following operator
 $A$ from the space
          $L^2(\B_n(\K),  dm(z))$
          to the space of functions
on $S^{rn-1}\times \R_+$:
\begin{equation}
Af(h,s)=\int_{\B_n(\K)}            f(z)
\frac{|1-\langle z,h \rangle|^{-(n+1)r/2+1+is}      }
     {|1-\langle z,z \rangle|^{(n+1)r/4+1/2+is/2} }   dz
.\end{equation}

\begin{lemma}
\label{l:intertwine}
The operator $A$ is intertwining, i.e.,
$$
A \rho(g)=\tau(g)A,\qquad \text{ for all $g$.}
$$
\end{lemma}

This statement is a useful two-step exercise.
First, it is worth to verify the lemma
in a straightforward way. Secondly,
it is interesting to find a way 
to invent the formula for the operator $A$
if you do not know it before.

\smallskip

{\bf\punct The Plancherel formula.}

 \begin{theorem}
\label{th:harmonic}
The operator $A$ is a unitary operator
\begin{equation}
L^2 (\B_n(\K),  dm(z))  \to L^2\left(S^{rn-1}\times \R_+,
 \left|\frac{\Gamma(b+is)\Gamma(c+is)}
{\Gamma(2is)}\right|^2 ds\,dh\right)
\label{eq:L2}
,\end{equation}
where
\begin{equation}
b=(n+1)r/4-1/2;\qquad c=(n-1)r/4+1/2
.\end{equation}
\end{theorem}

Keeping in the mind the previous lemma we get that
the operator
$A$ identifies the representation 
of the group  $\U(1,n;\K)$
in $L^2$ on the ball with a continuous direct sum
of representations of principal series.

\smallskip

{\sc Beginning of proof.} First, we explain the appearance
of the   $\Gamma$-factor. For this purpose we restrict
the operator 
 $A$  to the space of functions depending only
 on radius. It is convenient
 to define the variable
$$x=\frac{|h|^2}{1-|h|^2}$$
and set $f=f(x)$. Then the corresponding function
 $G(h,s)$ depends only the variable
 $s$, an uncomplicated calculation
 gives the familiar formula 
\begin{equation}
G(s)=\const\cdot\intt f(x) \FBC\densityx dx
.\end{equation}

Thus, we observe that the operator
 $A$ is a unitary operator from the space 
 of  $L^2$-functions on the ball depending only
 on radius to the space of functions on the half-cylinder depending only
 on $s$.

This is the main argument%
\footnote{For  $\K=\C$ there are also actions
of the group   $\U(1,n)$ given by
$$\rho(g)f(z)= f\bigl((a+zc)^{-1}(b+zd) \bigr)\,(a+zc)^k \ov{(a+zc)}^{\,\,-k} $$
 The problem of decomposition also is reduced to the index transform
(with another parameters), the representation has  finite 
discrete spectrum which is controlled by the Romanovsky
polynomials.}, it remains to apply some standard representation-theoretic
tricks.

\smallskip

{\bf\punct The end of the proof.} Denote
$G:=\U(1,n;\K)$, $K:=\U(1,\K)\times\U(n,\K)$. 
Denote the Hilbert spaces $L^2$ from row 
(\ref{eq:L2}) by  $V$ and $W$ respectively.
By $V^K$ and $W^K$ we denote the spaces
 of $K$-fixed vectors in $V$ and $W$.
By $P_V$ and $P_W$ we denote the projection operators 
to 
  $V^K$ and $W^K$.

Recall the following standard statement.

\begin{lemma}
  Let   $\rho(k)$ be a unitary representation of a compact group 
   $K$. Then the projection operator 
   to the space of $K$-fixed vectors 
   is given by the formula
$$P=\int_K\rho(k)\,dk,$$
where $dk$ is the Haar measure on  $K$,
normed in such a way that the measure of the whole
group is 1.
\end{lemma}

 \begin{corollary}
$P_W A=AP_V$.
\end{corollary}

\begin{lemma}
Any closed $G$-invariant subspace in $V$
contains a smooth function.
\end{lemma}

{\sc Proof.} Consider a sequence of smooth compactly supported positive functions
 $r_j$ on $G$ approximating the
$\delta$-function at unit.
For a vector $v\ne 0$ from the subspace
we get a sequence of smooth functions 
 $\int r_j(g)\rho(g)v\,dg$,
convergent to $v$. \hfill$\square$

\begin{lemma}
Any closed  $G$-invariant subspace in  $V$ contains a  
$K$-invariant vector.
\end{lemma}

{\sc Proof.} Consider a smooth function   $f$
from the subspace. Let $f(a)\ne 0$. Consider  $g\in G$
such that $0^{[g]}=a$. Next, average the function  $f(x^{[g]})$
by $K$. \hfill$\square$

\begin{corollary}
The linear span of vectors
 $\rho(g) v$, where  $g$ ranges in  $G$ and $v$
ranges in $V^K$, is dense in  $V$. 
\end{corollary}
 
{\sc Proof.} If not, we consider the orthogonal
comlement to this kinear span. It contains a 
 $K$-invariant vector. \hfill$\square$

\begin{lemma}   
\label{l:last}
The linear span of vectors
 $\tau(g) w$, where $g$ ranges in  $G$ and $w$
ranges in $W^K$, is dense in $W$. 
\end{lemma}

{\sc  Proof.} We use irreducibility
of representations 
$T_s$ of the principal series. As $w$ we take functions 
$f(x,s)=1$, if $|s-s_0|<\epsilon$ and 0 otherwise.
It is easy to show that there are no functions
orthogonal
to all
$\tau_g(g)f$.   \hfill$\square$

\smallskip

{\sc End of proof of Theorem \ref{th:harmonic}.}
Let $\rho(g)v$, $\rho(g')v'$ be as in the last lemma.
Then 
\begin{multline*}
\la \rho(g)v, \rho(g')v'\ra_V
\stackrel{ \text{(representation  $U$ is unitary)}}
=\la v, \rho(g^{-1}g')v'\ra_V \stackrel{\text{($P_V$ is projection operator)}} =
 \\
= \la v, P_V  \rho(g^{-1}g')v'\ra_V
\stackrel{ \text{(Plancherel formula)}}
=\la Av,   AP_W \rho(g^{-1}g')v'\ra_W=
\\
\stackrel{\text{($A$ is intertwining)}}
=\la Av, P_W \tau(g^{-1}g')A v'\ra_W
\stackrel{\text{($P_W$ is a projection operator)}}=
\\=
 \la Av,  \tau(g^{-1}g)A v'\ra_W
 \stackrel{\text{(representation $\tau$ is unitary)}}
 =  \la \tau(g)A v, \tau(g')A v'\ra_W
.\end{multline*}

Therefore $A$ is an isometry. On the other hand,
the image of  $A$ contains  $W^K$
(by Theorem   \ref{th:unitary}) and therefore contains the whole
$W$
(by Lemma (\ref{l:last})).

{\tt Math.Dept., University of Vienna,

 Nordbergstrasse, 15,
Vienna, Austria

\&

Institute for Theoretical and Experimental Physics,

Bolshaya Cheremushkinskaya, 25, Moscow 117259,
Russia

\&

Mech.Math. Dept., Moscow State University,
Vorob'evy Gory, Moscow

e-mail: neretin(at) mccme.ru

URL:www.mat.univie.ac.at/$\sim$neretin

wwwth.itep.ru/$\sim$neretin}

\end{document}